\documentclass[12pt]{article}
\usepackage[cp1251]{inputenc}
\usepackage[russian]{babel}
\usepackage[T2B]{fontenc}

\usepackage[pdftex, colorlinks, urlcolor=blue, pdfborder={0 0 0 [2 0]}]{hyperref} 

\usepackage{amsfonts,amsmath,amssymb}
\usepackage{amsthm}
\usepackage{array}
\usepackage{esint}
\usepackage{graphicx}

\usepackage{misccorr} 

\newtheorem{theorem}{\hspace{\parindent}\bf{Теорема}}

\newtheorem{m-lemma}{\hspace{\parindent}\bf{Микролемма}}

\newtheorem{dfn}{\hspace{\parindent}\sl{О\,п\,р\,е\,д\,е\,л\,е\,н\,и\,е\,}}

\newenvironment{proofs}
{\vspace{1pt}\par{\sl%
Д\,о\,к\,а\,з\,а\,т\,е\,л\,ь\,с\,т\,в\,о.\,\ }}%
{\noindent\vspace{1pt}}

{\vspace{1pt}\par{\sl%
Д\,о\,к\,а\,з\,а\,т\,е\,л\,ь\,с\,т\,в\,о\ \ т\,е\,о\,р\,е\,м\,ы\, }}%
{\noindent\vspace{1pt}}

{\noindent\vspace{1pt}}

\newenvironment{proof-m}
{\vspace{1pt}\par{\sl%
Д\,о\,к\,а\,з\,а\,т\,е\,л\,ь\,с\,т\,в\,о\ \ м\,и\,к\,р\,о\,л\,е\,м\,м\,ы.\,\ }}%
{\noindent\vspace{1pt}}

{\vspace{1pt}\par{\sl%
Р\,е\,ш\,е\,н\,и\,е.\,}}%
{\noindent\vspace{1pt}}

\makeatletter
\renewcommand{\@biblabel}[1]{#1.} 
\makeatother

\begin{document}

\title{Стохастические первые интегралы, ядра интегральных инвариантов и уравнения Колмогорова}

\author{В. А. Дубко, Е. В. Карачанская  \\
Академия муниципального управления, Киев\\
Тихоокеанский госуниверситет, Хабаровск}


\date{}

\maketitle
\begin{abstract}
В работе мы представляем стохастический первый интеграл, обобщенную формулу Ито--Вентцеля и еее применение для получения  уравнений для стохастического первого интеграла, ядер интегральных инвариантов и уравнений Колмогорова для плотности переходных вероятностей случайных процессов, описываемых обобщенным СДУ Ито.

\end{abstract}

{\it Ключевые слова}: стохастический первый интеграл, стохастическое ядро стохастического интегрального инварианта, локальное стохастическое ядро, обобщенное уравнение Ито, уравнения Колмогорова

\section*{Введение}

В аналитической механике понятие первого интеграла решения детерминированной динамической системы играет фундаментальную роль, поскольку он связан с законами сохранения (энергии, силы, массы и т.п.). Выяснилось, что  первый интеграл существует и в теории стохастических дифференциальных уравнений (СДУ). Это {\it первый интеграл} для системы СДУ Ито
(В. А. Дубко, 1978 \cite{D_78});
{\it первый прямой интеграл}  и {\it первый обратный интеграл} для системы СДУ Ито (Н.~В.~Крылов и Б.~Л.~Розовский, 1982 \cite{KR_82}), {\it стохастический первый интеграл} для системы обобщенных уравнений Ито (ОСДУ) \cite{D_02,11_KchIto}.
Однако классический перенос понятия первого интеграла на стохастические системы уравнений не возможен.

Источником инвариантов, первых интегралов являются различные законы сохранения (энергии, массы, импульса, момента импульса и т. п.). Например, если набор счетного числа начальных значение решений некоторого динамического уравнения связать с точками (аналогами частиц), то при выполнении условий существования и единственности решений число этих точек будет оставаться одним и тем же в любой момент времени. Предельным обобщением этого представления является плотность числа точек и, соответственно, неизменность интеграла от этой плотности по пространству. Функцию, обладающую таким свойством  называют ядром интегрального инварианта.

При определенных ограничениях возможно построить уравнение в частных производных для ядер \cite{D_02,D_89,13_KchArxDir}. Не исключается и ситуация, когда ограничениями являются свойства функционалов, сохраняющихся в некоторой ограниченной области координат и времени. Для таких случаев приходят к представлению о локальных инвариантах \cite{D_95}. Эволюционирующие структуры и функционалы, связанные с областью начальных значений, рассматриваются как {\it динамические инварианты.} Примерами динамических инвариантов могут служить элемент фазового объема, гиперповерхности.

Теория стохастических инвариантов представляет собой один из подходов к изучению стохастических динамических систем, описываемых уравнениями Ито \cite{D_89,D_84,D_98}. Этот подход оказался эффективен для установления наличия сохраняющихся функционалов у стохастических эволюционирующих структур (первых и стохастических первых интегралов уравнений Ито \cite{D_78,D_02-Dan}, длины случайной цепи \cite{10_KchTSP}, постоянной скорости случайно движущейся частицы \cite{97_DubChalBroun,04_ChConf}), построения точного решения стохастического дифференциального уравнения (СДУ) типа Ланжевена \cite{98_DubChalSolut}, определения моментных характеристик диффундирующей на сфере частицы \cite{10_KchHf,12_KchSphera}, получения  формулы Ито-Вентцеля и ее аналога -- обобщенной формулы Ито-Вентцеля для обобщенных СДУ Ито с винеровской и пуассоновской составляющей \cite{D_89,D_02,11_KchOboz,12_KchDubPrep}, а также  для обоснования возможности определения и построения программных управлений с вероятностью единица в стохастических системах, находящихся в условиях сильных возмущений \cite{07_ChOboz,09_ChUprW,11_KchUpr}.

Продемонстрируем, как можно использовать представление о стохастическом первом интеграле для построения и доказательства теорем о виде  уравнений для него, и теорем  существовании и единственности решений уравнений для ядер интегральных инвариантов. Важным моментом в том служат полученное в наших работах правило дифференцирования  случайных функций, зависящих от решений обобщенных уравнений Ито. Этому правилу мы дали название "<обобщенная формула Ито -- Вентцеля"> \cite{D_02,11_KchIto,12_KchDubPrep}. Выбор такого названия связан с тем, что при отсутствии пуассоновских возмущений полученная нами формула переходит в формулу, известную под названием "<формула Ито -- Вентцеля"> \cite{KR_82,D_89,Wentzel_65}.

Цель статьи: продемонстрировать возможности применения обобщенной  формулы Ито -- Вентцеля для построения уравнения для стохастического первого интеграла, строгого обоснования получения уравнения для стохастических ядер стохастических интегральных инвариантов и получения уравнений Колмогорова для переходных вероятностей \cite{13_DKchUch}.

Работа состоит из трех  частей. В первой мы определяем понятие стохастического первого интеграла, приводим обобщенную формулу Ито -- Вентцеля и уравнение для стохастического первого интеграла; во второй -- вводим понятие локальной стохастической плотности динамического инварианта, связанного с решением обобщенного уравнения Ито, устанавливаем вид уравнения для плотности и связь локальной стохастической плотности с понятием стохастического ядра стохастического интегрального инварианта;
в третьей --
строим уравнения Колмогорова.

\section{Стохастические первые интегралы}

Пусть $(\Omega,\mathcal{F}, \mathbf{P})$ -- полное вероятностное пространство, ${\rm {\mathcal F}}_{t}=\bigl\{{\mathcal F}_{t},  \ \ t\in [0,T]\bigr\}$,   ${\mathcal F}_{s}\subset {\mathcal F}_{t} $, $s<t$ -- неубывающий поток $\sigma-$алгебр.
Рассматриваются следующие случайные процессы: ${\bf w}(t)$ -- $m$-мерный винеровский процесс;
  $\nu (\Delta t;\Delta \gamma )$ -- стандартная мера Пуассона, определенная на $ [0;T]\times\mathbb{R}^{n'} $,  ${\rm M}[\nu (\Delta t,\Delta \gamma )]=\Delta t\cdot \Pi (\Delta \gamma )$; $\int_{{\mathbb R} (\gamma )}  \Pi (d\gamma )<\infty $, $ {\mathbb R} (\gamma )=\mathbb{R}^{n'}$ (это условие означает, что интенсивность скачков за бесконечно малый промежуток времени  конечна);
 одномерные винеровские процессы $w_{k}(t)$ и пуассоновская мера $\nu([0;T],\mathcal{A})$ определены на вероятностном пространстве $(\Omega,\mathcal{F},\mathbf{P})$, $\mathcal{F}_{t}-$измеримы для всех $t>0$ для любых множеств $\mathcal{A}$ из $\sigma-$алгебры борелевских множеств, и являются взаимно независимыми.

Пусть случайный процесс $x(t)$, определенный на ${\mathbb R}^{n} $, является решением системы  стохастических уравнений \cite{GS_68}:
\begin{equation} \label{Ayd01.1u}
\begin{array}{c}
\displaystyle{dx_{i} (t)=a_{i} (t )
dt+b_{i k} (t )
dw_{k} (t)+\int_{{\mathbb R} (\gamma )} g_{i} (t; \gamma ) \nu (dt;d\gamma ),} \\
\displaystyle {x(t)=x(t;x_{0})\bigr|_{t=0} =x_{0}, \ \ \ {\mbox {\rm для всех}} \ \ x_{0}\in {\rm {\mathbb R}}^{n} }, \ \ i=\overline{{1,n}};\ k=\overline{{1,m}},
\end{array}
\end{equation}
где
$i =\overline{{1,n}}$, $k=\overline{{1,m}} $, $ {\mathbb R} (\gamma )=\mathbb{R}^{n'}$ -- пространство, на котором определены векторы $\gamma$; и  по индексам, встречающимся дважды, ведется  суммирование. Коэффициенты $a(t;x)$, $b(t;x)$ и $g(t;x;\gamma )$ удовлетворяют условиям, которые  обеспечивают существование и единственности решения уравнения \eqref{Ayd01.1u} (см., например,  \cite{GS_68}).

Далее, следуя \cite{GS_68}, вместо  $\int _{{\rm {\mathbb R}}\left(\gamma \right)} $ будем использовать обозначение $\int  $.

Пусть $u(t;x;\omega )$  -- случайная функция, определенная на том же вероятностном пространстве, что и решение системы                         \eqref{Ayd01.1u}

\begin{dfn}{\rm \cite{D_02,11_KchIto}}\label{Ayddf1}
Случайную функцию $u(t; x ;\omega)$, определенную на
том же вероятностном пространстве, что и решение системы
\eqref{Ayd01.1u}, будем называть стохастическим
первым интегралом системы  ОСДУ Ито \eqref{Ayd01.1u}, если с
вероятностью единица выполняется условие
\begin{equation}\label{xGrindEQ__1_}
u\bigl(t; x (t;  x_{0});\omega\bigr)=u\bigl(0; x_{0}\bigr)
\end{equation}
для любого решения $ x (t; x_{0};\omega)$ системы
\eqref{Ayd01.1u}.
\end{dfn}

Символ $\omega$, как принято, опускаем.

Приведем теперь формулировку обобщенной формулы Ито -- Вентцеля, которая  будет использоваться нами далее, в том числе и для получения уравнения для стохастического первого интеграла.
\begin{theorem} {\rm(Обобщенная формула Ито -- Вентцеля)} {\rm \cite{13_KchArxDir,12_KchDubPrep}}\label{t1}
Пусть случайный процесс $ x (t)\in  {\mathbb R}^{n}$  подчинен системе \eqref{Ayd01.1u}
с коэффициентами, для которых  выполняются условия ${\bf \mathcal{L}_{1}}$:
\begin{description}
  \item[$ \mathcal{L}_{1}.1.$] $\displaystyle\int_{0}^{T}|a_{i} (t)| dt<\infty;$ \ \
 $\displaystyle\int _{0}^{T}|b_{i k} (t)|^{2} dt<\infty$;
  \item[$ \mathcal{L}_{1}.2.$] $\displaystyle\int_{0}^{T}dt\int   |g_{i} (t;\gamma)|^{s} \Pi (d\gamma)<\infty, \ \ s=1,2.$
\end{description}
Тогда, если $F(t; x )$, $(t; x )\in[0,T]\times\mathbb{R}^{n} $ -- скалярная функция, обобщенный стохастический дифференциал которой имеет вид:
\begin{equation} \label{GrindEQ__2_5_1_}
d_{t} F(t; x )=Q(t; x )dt+D_{k} (t; x )dw_{k} (t)+\int G(t; x ;\gamma ) \nu (dt;d\gamma )
\end{equation}
и для коэффициентов \eqref{GrindEQ__2_5_1_} выполнены условия ${\bf \mathcal{L}_{2}}$:
\begin{description}
  \item[$ \mathcal{L}_{2} .1.$] $Q(t; x )$, $D_{k}(t; x )$, $G(t; x ;\gamma)$ -- в общем случайные функции, измеримые  относительно  потока $\sigma-$алгебр $ \mathcal{F}_{t} $,
согласованного с процессами $ w_{k}(t)$, $k=\overline{1,m}$, и $\nu(t;\mathcal{A})$ из \eqref{Ayd01.1u} для любого множества  $\mathcal{A}\in \mathfrak{B}$  из фиксированной борелевской $\sigma-$алгебры  {\rm(\cite[с. 266]{GS_68})};
  \item[$ \mathcal{L}_{2}.2.$] $
 Q (t; x )\in \mathcal{C}_{t,x}^{1,2}, \ \ D_{k}(t; x )\in \mathcal{C}_{t,x}^{1,2}, \ \ \ G(t; x ;\gamma)\in \mathcal{C}_{t,x,\gamma}^{1,2,1}
$.
\end{description}
то существует стохастический дифференциал:
\begin{equation} \label{GrindEQ__2_5_2_}
\begin{array}{c}
 \displaystyle d_{t} F(t; x (t))=Q(t; x (t))dt+D_{k} (t; x (t))dw_{k} +\\
+\displaystyle \Bigl[a_{i} (t)\frac{\partial F(t; x )}{\partial x_{i} } +\frac{1}{2} b_{ik} (t)b_{j,k} (t)\frac{\partial ^{{\kern 1pt} 2} F(t; x )}{\partial x_{i} \partial x_{j} }\Bigr. + \\
+\Bigl.\displaystyle
b_{ik} (t)\frac{\partial D_{k} (t; x )}{\partial x_{i} }\Bigr]\bigr|_{ x = x (t)} dt+
b_{ik} (t)\displaystyle\frac{\partial F(t; x )}{\partial x_{i} }\bigr|_{ x = x (t)} dw_{k} + \\
+\displaystyle\int \Bigl[(F(t; x (t)+g(t;\gamma ))-F(t; x (t))\Bigr]\nu (dt;d\gamma ) +\\
+ \displaystyle\int G\bigl(t; x (t)+g(t;\gamma );\gamma\bigr)\nu (dt;d\gamma ). \end{array}
\end{equation}
\end{theorem}

Отметим, что обобщенную формулу Ито -- Вентцеля будем применять при  дополнительных ограничениях на коэффициенты уравнения \eqref{Ayd01.1u} \cite{11_KchOboz}, необходимых для дальнейших рассуждений:
\begin{equation}
 a_{i}(t)=a_{i} (t; x )\in \mathcal{C}_{t,x}^{1,2}, \ \ b_{k}(t)=b_{k}(t; x )\in \mathcal{C}_{t,x}^{1,2}, \ \ \ g_{i}(t;\gamma)=g_{i}(t; x ;\gamma)\in \mathcal{C}_{t,x,\gamma}^{1,2,1}.
\end{equation}

Для компактности записи далее будем использовать обозначение $\displaystyle \frac{\partial f(t;x(t))}{\partial x_{j} }$ вместо $\displaystyle \frac{\partial f(t;x )}{\partial x_{j} }\Bigl|_{x=x(t)}$.

Построим уравнение для стохастических первых интегралов обобщенного уравнения Ито, предполагая, что дифференциал стохастического первого интеграла существует и имеет вид:
\begin{equation} \label{xGrindEQ__2_}
d_{t} u(t;x)=Q(t;x)dt+D_{k} (t;x)dw_{k} (t)+\int G(t;x;\gamma ) \nu (dt;d\gamma )
\end{equation}

Опираясь на определение \ref{Ayddf1} и представление \eqref{xGrindEQ__2_}, применим обобщенную формулу Ито -- Вентцеля \eqref{GrindEQ__2_5_2_}:
\begin{equation} \label{xGrindEQ__3_}
\begin{array}{c}
  \displaystyle d_{t}u(t;x(t))=Q(t;x(t))dt+D_{k} (t;x(t))dw_{k} +b_{ik} (t)\frac{\partial }{\partial x_{i} } u(t;x(t))dw_{k} + \\
  \displaystyle  +\Bigl[a_{i} (t)\frac{\partial }{\partial x_{i} } u(t;x(t))+\frac{1}{2} b_{ik} (t)b_{j,k} (t)\frac{\partial ^{{\kern 1pt} 2} u(t;x(t))}{\partial x_{i} \partial x_{j} } + \\
  \displaystyle  +b_{ik} (t)\frac{\partial }{\partial x_{i} } D_{k} (t;x(t))\Bigr]dt+\int G(t;x(t)+g(t;\gamma );\gamma )\nu (dt;d\gamma )+\\
   \displaystyle  +\int \bigl[(u(t;x(t)+g(t;\gamma ))-u(t;x(t))\bigr]\nu (dt;d\gamma )=0.
\end{array}
\end{equation}

Вид уравнения для стохастического первого интеграла зависит от того, зависит или не зависит функция $g(t;\gamma )$ от $x$. Рассмотрим обе ситуации.

Пусть функция $g(t;\gamma )$ явно не зависит от  $x$, и коэффициенты уравнения \eqref{xGrindEQ__2_} имеют следующий вид
\begin{subequations}\label{xGrindEQ__4_}
\begin{align}
 \begin{array}{c}  \displaystyle{Q(t;x)=-\bigl[a_{i} (t)\frac{\partial }{\partial x_{i} } u(t;x)+\frac{1}{2} b_{ik} (t)b_{j,k} (t)\frac{\partial ^{ 2} u(t;x)}{\partial x_{i} \partial x_{j} } +} \\  \displaystyle{-b_{ik} (t)\frac{\partial }{\partial x_{i} } (b_{j,k} (t)\frac{\partial }{\partial x_{j} } u(t;x))\bigr]},
 \end{array}
 \label{xGrindEQ__4_a}\\
 \displaystyle D_{k} (t;x) = -b_{ik} (t)\frac{\partial }{\partial x_{i} } u(t;x),  \ \ \ \ \ \ \ \ \ \ \ \ \ \ \ \ \label{xGrindEQ__4_b} \\
   \displaystyle G(t;x;\gamma )=u(t;x-g(t;\gamma ))-u(x;t).\ \ \ \ \ \ \ \ \ \ \ \label{xGrindEQ__4_c}
\end{align}
\end{subequations}

\begin{theorem} Для того, чтобы непрерывная, ограниченная вместе со своими производными, вплоть до второй по компонентам $x$, функция $u(t,x)$, дифференциал от которой определяется \eqref{xGrindEQ__2_}, была первым стохастическим интегралом системы  \eqref{Ayd01.1u},  достаточно чтобы она являлась решением уравнения  \eqref{xGrindEQ__2_}, с коэффициентами  \eqref{xGrindEQ__4_}.
\end{theorem}

\begin{proofs} \underbar{Достаточность} связана с проверкой выполнения равенства \eqref{xGrindEQ__1_}. Это равенство заведомо будет реализовано, если коэффициенты уравнения \eqref{xGrindEQ__2_} будут определяться \eqref{xGrindEQ__4_}, что проверяется непосредственно  прямой подстановкой в \eqref{xGrindEQ__3_}.~$\lozenge$
\end{proofs}

Перейдем к рассмотрению случая, когда $g(t;\gamma )=g(t;x;\gamma )$:
\begin{equation} \label{xGrindEQ__5_}
dx(t)=a(t)dt+b_{k} (t)dw_{k} (t)+\int g(t;x(t);\gamma ) {\rm \; }\nu (dt;d\gamma ),
\end{equation}
где $a(t)=(a_{i}(t))$, $b_{k}(t)=(b_{ik}(t))$, $g(\cdot)=(g_{i}(\cdot))$, $i=\overline{{1,n}}$,  $k=\overline{{1,m}}$.

Введем новые функции и сделаем ряд замечаний относительно их свойств.

Обозначим через $x^{-1} (t;y;\gamma )$ решения относительно переменной $x$ равенства $y=x+g(t;x;\gamma )$  (см. \eqref{xGrindEQ__3_}).

Из этого определения следует, что если рассматривать области однозначности, то выбрав из такой области
\begin{equation}\label{6a}
 y=z+g(t;z;\gamma ),                                                        \end{equation}
приходим к тождеству:
\begin{equation}\label{6b}
x^{-1} (t;z+g(t;z;\gamma );\gamma )=z.                                                   \end{equation}

Пусть теперь $u(t;x)$ случайная функция, стохастический дифференциал для которой определяется \eqref{xGrindEQ__2_} с коэффициентами \eqref{xGrindEQ__4_a}, \eqref{xGrindEQ__4_b}  и
\begin{equation}\label{xGrindEQ__7_}
 G(t;x;\gamma )=u(t;x-g(t;x^{-1} (t;x;\gamma );\gamma ))-u(x;t).
\end{equation}

\begin{theorem} Случайная функция $u(t;x)\in \mathcal{C}_{t,x}^{1,2}$, стохастический дифференциал которой определяется \eqref{xGrindEQ__2_}, а коэффициенты $D_{k} (t;x)$, $Q(t;x)$ и $G(t;x;\gamma )$, соответственно,  \eqref{xGrindEQ__4_a}, \eqref{xGrindEQ__4_b}  и \eqref{xGrindEQ__7_}, является стохастическим первым интегралом обобщенного уравнения Ито \eqref{xGrindEQ__5_}. Причем  при указанных ограничениях на коэффицициенты эти условия являются необходимыми и достаточными.
\end{theorem}

\begin{proofs}  Для доказательства достаточно проверить равенство нулю суммы слагаемых в \eqref{xGrindEQ__3_}, связанных с интегралами по пуассоновской мере. Действительно, с учетом представления \eqref{xGrindEQ__7_} и свойств  \eqref{6a} и \eqref{6b}, убеждаемся, что:
\begin{equation}
\begin{array}{c}
   \displaystyle\int G(t;x(t)+g(t;x(t);\gamma ))\nu (dt;d\gamma )+ \\
  + \displaystyle \int \bigl[(u(t;x(t)+g(t;x(t);\gamma ))-u(t;x(t))\bigr]\nu (dt;d\gamma ) = \\
  = \displaystyle\int \bigl[(u(t;x(t)+g(t;x(t);\gamma )-g(t;x^{-1} (t;x+g(t;x(t);\gamma );\gamma )- \\
  -u(x(t)+g(t;x(t);\gamma );t)\bigr]\nu (dt;d\gamma )+\\
  + \displaystyle\int \bigl[(u(t;x(t)+g(t;x(t);\gamma ))-u(t;x(t))\bigr]\nu (dt;d\gamma ) \equiv 0.
\end{array}
\end{equation}

Достаточность доказана.

\underbar{Необходимость} следует из условий существования и единственности решений \eqref{xGrindEQ__5_} и, как следствие, единственности представления для стохастического дифференциала для $F(t;x(t))$ от произвольной случайной функции $F(t;x)$, стохастический дифференциал от которой определяется  \eqref{xGrindEQ__1_}. $\lozenge$
\end{proofs}

\section{Ядра инвариантов для обобщенных уравнений Ито и уравнения для них}

Пусть теперь случайный процесс $x(t)$, определенный на ${\mathbb R}^{n} $, является решением системы стохастических дифференциальных уравнений  при более строгих ограничениях на часть коэффициентов \cite[с.~278--290, 298--302]{GS_68}:
\begin{equation} \label{Ayd01.1}
\begin{array}{c}
\displaystyle{dx_{i} (t)=a_{i} (t;x(t))
dt+b_{ik} (t;x(t))
dw_{k} (t)+\int g_{i} (t;x(t);\gamma ) \nu (dt;d\gamma ),} \\
\displaystyle {x(t)=x(t;x_{0})\bigr|_{t=0} =x_{0}, \ \ \ {\mbox {\rm для всех}} \ \ x_{0}\in {\rm {\mathbb R}}^{n} }, \ \ i =\overline{{1,n}}, \ \ k=\overline{{1,m}},
\end{array}
\end{equation}
 \begin{equation}\label{exiseq}
 a_{i}(t; x )\in \mathcal{C}_{t,x}^{1,1}, \ \ b_{ik}(t; x )\in \mathcal{C}_{t,x}^{1,2};
\end{equation}
и
\begin{equation} \label{pu}
\begin{array}{c}
\displaystyle{\int _{0}^{T}dt \int  |\nabla ^{\alpha } g(t;x;\gamma )| ^{\beta } \Pi (d\gamma )<\infty ,\ \ \ \alpha = 1,2 , \ \ \beta =\overline{1,4}, }
 \end{array}
\end{equation}
где $\nabla ^{k} $ обозначает всевозможные комбинации частных производных $k$-го порядка по компонентам $x$  в предположении, что соответствующие выражения непрерывны по совокупности переменных.

Условий  \eqref{exiseq} и \eqref{pu} достаточно для существования и единственности решения уравнения \eqref{Ayd01.1} \cite{GS_68}.

Введем несколько определений.

\begin{dfn}{\rm \cite{D-dis_79}}
Пусть $S(t)=S(t;v)$, $v\subset\hat{\mathcal{F}}\subset\mathcal{F}$ -- измеримое отображение $\mathbb{R}^{n}$ в $\mathbb{R}^{n}$. Множество
$
\Delta (t) = S(t;v)\cdot \Delta(0)
$
будем называть динамическим инвариантом области $\Delta(0)$ для процесса $x(t)\in \mathbb{R}^{n}$, если
$$
\mathbf{P}\Bigl\{ x(t)\in \Delta(t)\bigl| x=x_{0}\Bigr\}=1, \ \ \ {\mbox {\rm для всех}} \ \ t>0, \ \ \ {\mbox {\rm для всех}} \ \ x\in\Delta(0).
$$
\end{dfn}

Пусть $\rho (t;x;\omega )$ -- случайная функция, измеримая относительно потока   $\sigma$-алгебр  ${\rm {\mathcal F}}_{t}$, согласованного с процессами ${w}_{k}(t)$, $k=\overline{1,m}$ и $\nu (t;\Delta \gamma )$ (в дальнейшем индекс $\omega $ будем опускать).

\begin{dfn}
Случайную функцию  $\rho (t;x)$  назовем локальной стохастической плотностью динамического инварианта для уравнения, связанного с обобщенным уравнением Ито \eqref{GrindEQ_D_2_}, если случайная функция $u(t;J;x)=J\rho (t;x)$ является стохастическим первым интегралом, т. е. удовлетворяет равенству
\begin{equation}\label{GrindEQ_D_2_}
J(t;x_{0})\rho (t;x(t;x_{0}))=\rho (0;x_{0}),\ \ \ {\mbox {\rm для всех}}\ \ x_{0}\in \Gamma\subset R^{n},
\end{equation}
где функция $J (t )=J (t,x_{0} )$ является решением уравнения
\begin{equation} \label{GrindEQ_D_3_}
\begin{array}{c}
\displaystyle{dJ (t )=J (t ) \left\{{\rm K}(t)dt+\frac{\partial b_{ik} (t;x(t) )}{\partial x_{i}} dw_{k} (t ) \right.+}\\
 +\displaystyle\left.\int \Bigl(\det \bigl[A(\delta _{i,j} +\frac{\partial g_{i} (t;x (t);\gamma )}{\partial x_{j} } )\bigr]-1\Bigr)\nu (dt,d\gamma ) \right\}, \\
\displaystyle{J (t )=J (t,x_{0} )\bigl|_{t=0} =1,{\rm \; \; dim}A(\cdot )=n\times n;}
\end{array}
\end{equation}
$\delta _{i,j} $ -- символ Кронекера;
$$
\begin{array}{c}
K(t)=\displaystyle\Bigl[\frac{\partial a_{i}(t;x(t) )}{\partial x_{i} } + \\
\displaystyle  +\frac{1}{2} \left(\frac{\partial b_{ik} (t;x(t) )}{\partial x_{i} } \cdot \frac{\partial b_{j,k} (t;x (t) )}{\partial x_{j} } -\frac{\partial b_{ik} (t;x(t) )}{\partial x_{j} } \cdot \frac{\partial b_{ik} (t;x(t) )}{\partial x_{i} } \right) .
\end{array}
$$
а $x(t)$ -- решение уравнения \eqref{Ayd01.1}.
\end{dfn}

Уравнение \eqref{GrindEQ_D_3_} является уравнением для якобиана преобразования от $ x_{0}$ к $ x(t;x_{0})$. Решение \eqref{GrindEQ_D_3_} может быть представлено в виде:
$$
\begin{array}{c}
J(t)=\displaystyle\exp \Bigl\{\int_{0}^{t}\Bigl[K(\tau ) -\frac{1}{2} \left(\frac{\partial b_{ik} (\tau ;x(\tau) )}{\partial x_{i} } \right)^{2} \Bigr]d\tau +\int_{0}^{t}\frac{\partial b_{ik}  (\tau ;x (\tau) )}{\partial x_{i} } dw_{k} (\tau ) +    \\
   +\displaystyle \int _{0}^{t}\int
   \ln \Bigl|\det \bigl[A(\delta _{i,j} +\frac{\partial g_{i} (\tau ;x(\tau)
   ;\gamma )}{\partial x_{j} } )\bigr]\Bigr|\nu (d\tau ,d\gamma ) .
\end{array}
$$

Отметим, что $J (t )>0$ для всех $t\ge 0$.

Перейдем к установлению вида уравнения в частных производных для функции $\rho (t;x)$.
Из равенства \eqref{GrindEQ_D_2_} следует, что
\begin{equation} \label{GrindEQ_D_4_}
d_{t}J(t;x_{0})\rho (t;x(t;x_{0}))=0, \ \ \ {\mbox {\rm для всех}} \ \ x_{0}\in \Gamma \subset  {\mathbb R} ^{n}
\end{equation}

Пусть функция $\rho (t;x)$ является решением уравнения
\begin{equation} \label{GrindEQ_D_5_}
\begin{array}{c}
\displaystyle{d_{t} \rho (t;x)=Q(t;x)dt+D_{k} (t;x)dw_{k} (t)+\int _{}G(t;x;\gamma ) \nu (dt;d\gamma ),} \\
\displaystyle{\rho (t;x)\left|{}_{t=0} =\right. \rho (x)\in C_{0}^{2} .}
\end{array}
\end{equation}
Здесь роль начальных условий выполняет функция. Относительно коэффициентов этого уравнения предполагаем, что
они удовлетворяют условиям
\begin{equation}\label{exiseqQQ}
 Q(t; x )\in \mathcal{C}_{t,x}^{1,1}, \ \ D_{k}(t; x )\in \mathcal{C}_{t,x}^{1,2};
\end{equation}
и
\begin{equation} \label{puGG}
\begin{array}{c}
\displaystyle{\int _{0}^{T}dt \int|\nabla ^{\alpha } G(t;x;\gamma )| ^{\beta } \Pi (d\gamma )<\infty ,\ \ \ \alpha = 1,2 , \ \ \beta =\overline{1,4}}
 \end{array}
\end{equation}

Покажем, что $\rho (t;x)$ должно являться решением стохастического уравнения в частных производных следующего вида, полученное нами ранее другими методами \cite{12_KchDubPrep}:
\begin{equation} \label{GrindEQ_D_6_}
\begin{array}{c}
\displaystyle{d _{t} \rho (t;x)=-\left[\frac{\partial \rho (t;x)a_{i} (t;x)}{\partial x_{i} } -\frac{\partial ^{2} \rho (t;x)b_{ik} (t;x)b_{j,k} (t;x)}{\partial x_{i} \partial x_{j} } \right]dt-} \\
\displaystyle{-\frac{\partial \rho (t;x)b_{ik} (t;x)}{\partial x_{i} } dw_{k} (t)+} \\
\displaystyle{+\int \Bigl[\rho \bigl(t;x-g(t;x^{-1} (t;x;\gamma );\gamma )\bigr)\bar{D}\bigl(x^{-1} (t;x;\gamma )\bigr)-\rho (t;x)\Bigr] \nu (dt;d\gamma ),} \\
\displaystyle{\rho (t;x)\bigl|_{t=0} =  \rho (x)\in C_{0}^{2} ; \ \ \ \rho (x)\ge 0, }
\end{array}
\end{equation}
где $x^{-1} (t;x;\gamma )$  определяется как решение относительно переменной  $y$ уравнения
\begin{equation}\label{GrindEQ_D_7_}
y+g(t;y;\gamma )=x
\end{equation}
во всех подобластях однозначности, а $\bar{D}\left(x^{-1} (t;x;\gamma )\right)$  -- якобиан перехода, соответствующий такой замене. Черточкой над большими буквами и далее будем обозначать детерминанты матриц.

Поскольку стохастический дифференциал -- обобщенное СДУ Ито -- можно представить в виде суммы слагаемых, одно из которых определяется наличием винеровского процесса, а второе -- наличием скачков пуассоновского процесса, то можно заменить требование \eqref{GrindEQ_D_4_} такими:
\begin{equation} \label{GrindEQ_D_8_}
\bigl[d_{t}J(t;x_{0})\rho \bigl(t;x(t;x_{0})\bigr)\bigr]_{1} =0,\ \ \ {\mbox {\rm для всех}} \ \  x_{0}\in \Gamma,
\end{equation}
\begin{equation} \label{GrindEQ_D_9_}
\bigl[d_{t}J(t;x_{0})\rho \bigl(t;x(t;x_{0})\bigr)\bigr]_{2} =0, \ \ \ {\mbox {\rm для всех}} \ \  x_{0}\in \Gamma ,
\end{equation}
где \eqref{GrindEQ_D_8_} соответствует стохастическому дифференциалу Ито, а \eqref{GrindEQ_D_9_}  связано с учетом пуассоновских компонент.

\begin{theorem} Пусть коэффициенты уравнений \eqref{Ayd01.1} и \eqref{GrindEQ_D_5_} удовлетворяют условиям \eqref{exiseq}, \eqref{pu} и  \eqref{exiseqQQ}, \eqref{puGG}  соответственно. Тогда коэффициенты $Q(t;x)$, $D_{k} (t;x)$, $G(t;x;\gamma )$ обеспечивающие выполнение условия \eqref{GrindEQ_D_4_}, однозначно определяются равенствами:
\begin{description}
  \item[$I.1.$] $-Q(t;x)=\displaystyle \frac{\partial \rho (t;x)a_{i} (t;x)}{\partial x_{i} } -\frac{\partial ^{2} \rho (t;x)b_{ik} (t;x)b_{j,k} (t;x)}{\partial x_{i} \partial x_{j} };$
  \item[$I.2.$] $-D_{k} (t;x)=\displaystyle \frac{\partial \rho (t;x)b_{ik} (t;x)}{\partial x_{i} };$
  \item[$I.3.$] $G(t;x;\gamma )= \displaystyle \rho \left(t;x-g(t;x^{-1} (t;x;\gamma );\gamma )\right)\bar{D}\left(x^{-1} (t;x;\gamma )\right)-\rho (t;x).$
\end{description}
Кроме того, при начальных условиях для \eqref{GrindEQ_D_6_} функция  $\rho (t;x)$ является единственным решением стохастического уравнения \eqref{GrindEQ_D_5_} с коэффициентами, определяемыми условиями I.1, I.2.

\end{theorem}

\begin{proofs} Ограничения теоремы на коэффициенты уравнения \eqref{GrindEQ_D_5_} обеспечивают применимость обобщенной формулы Ито -- Вентцеля \eqref{GrindEQ__2_5_2_}. В результате ее применения получим:
\begin{equation} \label{durad}
\begin{array}{c}
  \displaystyle d_{t}\rho (t;x(t))=Q(t;x(t))dt+D_{k} (t;x(t))dw_{k} +b_{ik} (t;x(t))\frac{\partial }{\partial x_{i} } \rho (t;x(t))dw_{k} + \\
  \displaystyle {+\Bigl[a_{i} (t;x(t))\frac{\partial }{\partial x_{i} } \rho (t;x(t))+\frac{1}{2} b_{ik} (t;x(t))b_{j,k} (t;x(t))\frac{\partial ^{ 2} \rho (t;x(t))}{\partial x_{i} \partial x_{j} } +} \\
 \displaystyle {+b_{ik} (t;x(t))\frac{\partial }{\partial x_{i} } D_{k} (t;x(t))\Bigr]dt+\int G\bigl(t;x(t)+g(t;x(t);\gamma )\bigr)\nu (dt;d\gamma )+}  \\
  +\displaystyle \int \Bigl[\bigl(\rho (t;x(t)+g(t;x(t);\gamma)\bigr)-\rho (t;x(t))\Bigr]\nu (dt;d\gamma).
\end{array}
\end{equation}
Рассмотрим составляющую дифференциала Ито этого выражения:
\begin{equation} \label{GrindEQ_D_10_}
\begin{array}{c}
 \bigl[d_{t}\rho (t;x(t))\bigr]_{1} =Q(t;x(t))dt+D_{k} (t;x(t))dw_{k}(t) +\\
  \displaystyle {+\bigl[a_{i} (t;x(t))\frac{\partial }{\partial x_{i} } \rho (t;x(t))+\frac{1}{2} b_{ik} (t;x(t))b_{j,k} (t;x(t))\frac{\partial ^{{\kern 1pt} 2} \rho (t;x(t))}{\partial x_{i} \partial x_{j} } +} \\
  \displaystyle {+b_{ik} (t;x(t))\frac{\partial }{\partial x_{i} } D_{k} (t;x(t))\bigr]dt+b_{ik} (t;x(t))\frac{\partial }{\partial x_{i} } \rho (t;x(t))dw_{k}(t) .}
\end{array}
\end{equation}

На явном виде $\bigl[d_{t}\rho (t;x(t))\bigr]_{2} $ мы остановимся чуть позднее.

Воспользовавшись \eqref{GrindEQ_D_3_}, \eqref{GrindEQ_D_10_} и формулой Ито \cite{GS_68}, находим:
$$
\begin{array}{c}
  \bigl[d_{t}J(t)\rho (t;x(t))\bigr]_{1} =\rho (t;x(t))\bigl[dJ(t)\bigr]_{1} +J(t)\bigl[d\rho (t;x(t))\bigr]_{1} + \\
  \displaystyle   +J(t)b_{ik} (t;x(t))\frac{\partial b_{j,k} (t;x(t))}{\partial x_{j} } \frac{\partial \rho (t;x(t))}{\partial x_{i} } dt= \\
  =J(t)  \displaystyle \Bigl[b_{ik} (t;x(t))\frac{\partial b_{j,k}  (t;x(t))}{\partial x_{j} } \frac{\partial \rho (t;x(t))}{\partial x_{i} } +a_{i} (t;x(t))\frac{\partial }{\partial x_{i} } \rho (t;x(t))+ \\
    \displaystyle {+\frac{1}{2} b_{ik} (t;x(t))b_{j,k} (t;x(t))\frac{\partial ^{2} \rho (t;x(t))}{\partial x_{i} \partial x_{j} } +b_{ik} (t;x(t))\frac{\partial }{\partial x_{i} } D_{k} (t;x(t))+} \\
    \displaystyle   {+b_{ik} (t;x(t))\frac{\partial }{\partial x_{i} } D_{k} (t;x(t))+Q(t;x(t))+}\\
  +  \displaystyle \rho (t;x(t))\frac{\partial a_{i} (t;x(t))}{\partial x_{i} } +\frac{1}{2} \bigl[\frac{\partial b_{ik}  (t;x(t))}{\partial x_{i} }  \frac{\partial b_{j,k}  (t;x(t))}{\partial x_{j} } -\frac{\partial b_{ik}  (t;x(t))}{\partial x_{j} }   \frac{\partial b_{ik}  (t;x(t))}{\partial x_{i} } \bigr]\Bigr]dt+\\
  +J(t)  \displaystyle \Bigl[b_{ik} (t;x(t))\frac{\partial \rho (t;x(t))}{\partial x_{i} } +\rho (t;x(t))\frac{\partial b_{ik} \left(t;x(t)\right)}{\partial x_{i} } +D_{k} (t;x(t))\Bigr]dw_{k} (t).
  \end{array}
$$
Учитывая \eqref{GrindEQ_D_8_} -- требование  обращения в ноль этого выражения, сгруппировав слагаемые, приходим к установлению необходимости выполнения  равенств I.1, I.2 теоремы.

Решения уравнений  \eqref{GrindEQ_D_5_}, \eqref{GrindEQ_D_6_} и \eqref{GrindEQ_D_10_} при указанных ограничениях, существуют и единственны \cite{KR_82}. Последнее и приводит к утверждению теоремы без пуассоновской составляющей.

Перейдем к рассмотрению условия I.3 теоремы.

На всех участках между пуассоновскими скачками, процесс подвержен только винеровским возмущениям и сохраняется порядок гладкости функции $\rho (t;x(t))$. В момент скачка происходит добавление функционала
\begin{equation}\label{Idt}
\begin{array}{c}
  \mathcal{I}( t)= \displaystyle\int_{0}^{t} \int G\bigl(\tau;x(\tau)+g(\tau;x(\tau);\gamma )\bigr)\nu (d\tau;d\gamma )+  \\
  +\displaystyle \int_{0}^{t}\int \Bigl[\bigl(\rho (\tau;x(\tau)+g(t;x(t);\gamma)\bigr)-\rho (\tau;x(\tau))\Bigr]\nu (d\tau;d\gamma)
 \end{array}
\end{equation}
неизменяющего  условия  гладкости коэффициентов.

В силу теорем о свойствах решений уравнений в частных производных \cite{KR_82}, и на последующем интервале между скачками решение будет существовать и не нарушится порядок его гладкости.

Поэтому остается проверить, что для системы
\begin{equation}\label{dI2}
\begin{array}{c}
\displaystyle {\bigl[d_{t} \rho (t;x(t))\bigr]_{2} =\int \Bigl\{ (\rho (t;x(t)+g(t;x(t);\gamma ))-\rho (t;x(t))+ } \\
\displaystyle {+\Bigl[\rho (t;x(t)+g(t;x(t);\gamma )-g(t;x^{-1} (t;x(t)+g(t;x(t);\gamma );\gamma );\gamma ))\Bigr]\times } \\
\displaystyle {\times \bar{D}\bigl(x^{-1} (t;x(t)+g(t;x(t);\gamma ))\bigr)-\rho (t;x(t)+g(t;x(t);\gamma ))\Bigr\} \nu (dt;d\gamma )=}\\
 =-\displaystyle \int \Bigr[\rho (t;x(t))-\rho (t;x(t))\bar{D}\bigl(\frac{\partial x_{i}^{-1} (t;x+g(t;x;\gamma );\gamma )}{\partial (x_{j} +g_{j} (t;x;\gamma ))} \bigr)\Bigl] \nu (dt;d\gamma ),\\
\displaystyle  \bigl[dJ(t)\bigr]_{2} =J(t)\int \Bigl(\det \Bigl[A(\delta _{i,j} +\frac{\partial g_{i} (t;x(t);\gamma )}{\partial x_{j} } )\Bigr]-1\Bigr)\nu (dt,d\gamma )
\end{array}
\end{equation}
будет обращаться в ноль в силу условия \eqref{GrindEQ_D_9_}, дифференциал:
$$
\begin{array}{c}
  \bigl[d_{t} J(t)\rho (t;x(t))\bigr]_{2} = \\
  =\displaystyle  \int \Bigl\{ \bigr[\rho (t;x(t))+(\rho (t;x(t))\bar{D}\bigl(\frac{\partial x_{j}^{-1} (t;x+g(t;x;\gamma );\gamma )}{\partial (x_{i} +g_{i} (t;x;\gamma ))} \bigr)-\rho (t;x(t))\bigr] \times \\
  \times \displaystyle \bigl [J(t)+J(t)\bar{A}\bigl(\delta _{i,j} (t)+\frac{\partial g_{i} (t;x(t);\gamma )}{\partial x_{l} } \bigr)-J(t))\bigr])-J(t)\rho (t;x(t))\Bigr\} \nu (dt;d\gamma )= \\
  =J(t)\displaystyle  \int \Bigl[\rho (t;x(t)) \bar{D}\bigl(\frac{\partial x_{j}^{-1} (t;x(t)+g(t;x;\gamma );\gamma )}{\partial (x_{i} +g_{i} (t;x;\gamma ))} \bigr)\bar{A}\bigl(\frac{\partial (g_{i} (t;x(t);\gamma )+x_{i} )}{\partial x_{j} } \bigr)-\\
  - \displaystyle  \rho (t;x(t))\Bigr]\nu (dt;d\gamma ).
\end{array}
$$

Подынтегральное выражение в первом из уравнений системы \eqref{dI2} получено сопоставлением уравнений \eqref{aGrindEQ__05_} и  \eqref{aGrindEQ__06_} для ядра и подстановкой в функционал \eqref{Idt} соответствующего выражения с учетом условия \eqref{GrindEQ_D_7_}.

Для завершения проверки I.3 достаточно убедиться в равенстве нулю подынтегральной разности. Это будет верно, если произведение детерминантов $\bar{D}(\cdot )\bar{A}(\cdot )=1$.

Учитывая, следующее из определения \eqref{GrindEQ_D_7_} то, что $x^{-1} \bigl(t;x+g(t;x;\gamma )\bigr)=x$, находим:
$$
\begin{array}{c}
  \displaystyle \det \Bigl[D\bigl(\frac{\partial x_{i}^{-1} (t;x+g(t;x;\gamma );\gamma )}{\partial (x_{l} +g_{i} (t;x;\gamma ))} \bigr)A \bigl(\frac{\partial (g_{l} (t;x;\gamma )+x_{i} )}{\partial x_{j} } \bigr)\Bigr]= \\
  =\displaystyle  \det   S\bigl(\frac{\partial x_{i}^{-1} (t;x+g(t;x;\gamma );\gamma )}{\partial (x_{l} +g_{l} (t;x;\gamma ))} \frac{\partial (g_{l} (t;x;\gamma )+x_{l} )}{\partial x_{j} } \bigr)= \\
  =\displaystyle \det  S\bigl(\frac{\partial x_{i}^{-1} (t;x+g(t;x;\gamma );\gamma )}{\partial x_{j} } \bigr)=\det S\bigl(\delta _{i,j} \bigr)=1.
\end{array}
$$

Это равенство и приводит к утверждению, что в области, где обеспечивается однозначное соответствие между переменными уравнения $y=x+g(t;x;\gamma )$, решение уравнения с коэффициентами, определяемыми равенствами теоремы, существует и единственно. $\lozenge$
\end{proofs}

При требовании существования и единственности решения уравнения \eqref{Ayd01.1} во всем пространстве, добавляются глобальные равенства \cite{D_02,D_89}:
\begin{equation}\label{ya}
\displaystyle \int _{{\rm {\mathbb R}}^{n} }f(x)\rho (t;x)d \Gamma (x)=\int _{{\rm {\mathbb R}}^{n} }f(x(t;y))\rho (y)d \Gamma(y), \ \ \ \int _{{\rm {\mathbb R}}^{n} }\rho (t;x)d \Gamma (x)=1
\end{equation}
для любой непрерывной и ограниченной $f(x)$, где $\rho (x)= \rho (0;x)$ удовлетворяет требованиям:
$$
\displaystyle \lim\limits_{|x|\to \infty} \frac{\partial^{k} \rho (t;x)}{\partial x_{i}^{k}}\Bigl|_{t=0}=0, \ \ \ \ k=0,1,2. \ \ \  i=\overline{1,n}.
$$
 Соотношения \eqref{ya} есть определение стохастического ядра стохастического интегрального инварианта \cite{D_02,12_KchDubPrep}.

Отметим, что как и для детерминированных уравнений существует полная система ядер интегральных инвариантов $\rho _{r} (t;x)$, $r=\overline{1,n+1}$, т. е. такая совокупность ядер, что  любая другая функция,
являющаяся ядром интегрального инварианта $n$-го порядка, может
быть представлена как функция от элементов этой совокупности.
Любое ядро $\rho (t;x)$  однозначно выражается через начальное значение $\rho (x)$ и набор этих ядер \cite{11_KchIto,D_89,Zubov_82}. Для стохастических процессов схема этого доказательства этого утверждения вполне аналогична доказательству для случая детерминированных систем \cite{Zubov_82}.

Подчеркнем, что существование плотности $\rho (t;x)$ связано со свойствами меры $\mu (t ;\Delta)$ ($\int _{{\rm {\mathbb R}}^{n} }\mu (t ;d(\Delta)) =const$), индуцируемой некоторым динамическим отображением. При условии, что существует, в каком-то смысле, предел $\mu (t ;\Delta)/\mu (\Delta )$ по $\Delta \downarrow 0$, он отождествляется с функцией $\rho (t;x)$. Но функция $\rho (t;x)$ может и не обладать необходимой гладкостью, позволяющей построить уравнение вида \eqref{GrindEQ_D_6_}.

\section{Применение уравнения для стохастической плотности обобщенного уравнения Ито  для построения уравнений Колмогорова}

Пусть $x(t)$ -- решение системы обобщенных уравнений Ито \eqref{Ayd01.1}, коэффициенты  $a_{i} (t;x)$,\, $b_{ik} (t;x)$  и $g(t;x;\gamma )$ которого удовлетворяют условиям \eqref{exiseq} и \eqref{pu}. Этих условий достаточно для существования и единственности решений уравнения \eqref{Ayd01.1} при некоторых дополнительных ограничениях и для рассматриваемых ниже решений уравнений.

Пусть неотрицательная случайная функция $\rho (t;x)$, определенная на том же потоке $\sigma -$алгебр, что и процесс $x(t)$, является стохастическим ядром интегрального инварианта, согласованного с процессом $x(t)$, т. е. удовлетворяет условиям \eqref{ya}.

Уравнения для $\rho (t;x)$ можно построить разными способами и для разных представлений стохастических уравнений. Мы воспользуемся уравнением для стохастических ядер со случайной мерой Пуассона \eqref{GrindEQ_D_6_}.

Перейдем в уравнении для ядер \eqref{GrindEQ_D_6_} к центрированной случайной мере $\tilde{\nu }(\Delta t,\Delta \gamma )=\nu (\Delta t,\Delta \gamma )-\Delta t\Pi (\Delta \gamma )$:
\begin{equation} \label{aGrindEQ__04_}
\begin{array}{c}
 \displaystyle d_{t} \rho (t;x)=-\frac{\partial \rho (t;x)b_{ik} (t;x)}{\partial x_{i} } dw_{k} (t)-\Bigl[\frac{\partial \left(\rho (t;x)a_{i} (t;x)\right)}{\partial x_{i} } - \\
  -\displaystyle\frac{1}{2} \frac{\partial ^{  2} \left(\rho (t;x)b_{ik} (t;x)b_{j,k} (t;x)\right)}{\partial x_{i} \partial x_{j} } \Bigr]dt+  \\
  +\displaystyle\int \bigl[\rho \left(t;x-g(t;x^{-1} (t;x;\gamma );\gamma )\right)  \bar{D} \left(x^{-1} (t;x;\gamma )\right)-\rho (t;x)\bigr]\Pi (d\gamma )dt+\\
  +\displaystyle\int \bigl[\rho \left(t;x-g(t;x^{-1} (t;x;\gamma );\gamma )\right)  \bar{D} \left(x^{-1} (t;x;\gamma )\right)-\rho (t;x)\bigr]\tilde{\nu }(dt;d\gamma ).
\end{array}
\end{equation}

Введем обозначение: $  \mathbf{M} [\rho (t;x)]=p(t;x)$.
Вычислив математическое ожидание от обеих частей в \eqref{aGrindEQ__04_}  находим, что $p(t;x)$ удовлетворяет уравнению
\begin{equation} \label{aGrindEQ__05_}
\begin{array}{c}
\displaystyle \frac{\partial p(t;x)}{\partial t} =-\frac{\partial \left(p(t;x)a_{i} (t;x)\right)}{\partial x_{i} } +\frac{1}{2} \frac{\partial ^{  2} \left(p(t;x)b_{ik} (t;x)b_{j,k} (t;x)\right)}{\partial x_{i} \partial x_{j} } +  \\
+\displaystyle
 \int [p\left(t;x-g(t;x^{-1} (t;x;\gamma );\gamma )\right) \bar{D} \left(x^{-1} (t;x;\gamma )\right)-p(t;x)]\Pi (d\gamma ).                 \end{array}
\end{equation}

Соотношение \eqref{aGrindEQ__05_} представляет собой  {\it уравнение Колмогорова для плотности}.

Построим уравнения для плотности переходных вероятностей процессов, являющихся решением обобщенных уравнений Ито \eqref{Ayd01.1}.

Под плотностью переходных вероятностей для случайного процесса понимают детерминированную функцию $p(t;x/s;y)$, обеспечивающую выполнение интегрального равенства
\begin{equation} \label{aGrindEQ__06_}
\displaystyle p(t;x)=\int _{{\mathbb R}^{n} } p(t;x/s;y)p(s;y) d\Gamma (y),\ \ \ t>s.
\end{equation}

Из уравнения \eqref{aGrindEQ__06_} следует, что функция $p(t;x/s;y)$ не зависит от распределения $p(s;y)$ для всех  $s\ge 0$, следовательно, и от произвольного $p(0;y)=\rho (y)$.

После подстановки \eqref{aGrindEQ__06_} в уравнение \eqref{aGrindEQ__05_} -- в силу произвольности $p(s;y)$ --  приходим к выводу, что требование \eqref{aGrindEQ__06_} будет выполнено, если
\begin{equation*}
\begin{array}{c}
\displaystyle \frac{\partial p(t;x/s;y)}{\partial t} = \ \ \ \\
=-\displaystyle \Bigl[\frac{\partial \left(p(t;x/s;y)a_{i} (t;x)\right)}{\partial x_{i} } -\frac{1}{2} \frac{\partial ^{  2} \left(p(t;x/s;y)b_{ik} (t;x)b_{j,k} (t;x)\right)}{\partial x_{i} \partial x_{j} } \Bigr]+
\\
\displaystyle {+\int \Bigl[p\left(t;x-g(t;x^{-1} (t;x;\gamma );\gamma )/s;y\right)  \bar{D} \left(x^{-1} (t;x;\gamma )\right)-} \\ -p(t;x/s;y)\Bigr]\Pi (d\gamma ).
\end{array}
\end{equation*}

Это -- {\it прямое уравнение Колмогорова} для плотности переходных вероятностей, связанной с уравнением \eqref{Ayd01.1}.

Перейдем к построению обратного уравнения для $p(t;x/s;y)$.

 Учитывая, что $p(t;x/s;y)$ -- детерминированная функция, можно представить \eqref{aGrindEQ__06_} в таком виде:
$$
p(t;x)=\int _{{\mathbb R}^{n} } \mathbf{M}[ p(t;x/s;y)\rho (s;y)]d\Gamma(y).
$$

От этого представления  с учетом свойства \eqref{ya} функции $\rho (s;y)$ переходим к равенству
$$
\begin{array}{c}
\displaystyle{p(t;x)=\int _{{\mathbb R}^{n} } p(t;x/s;y) \mathbf{M}[\rho (s;y)]d\Gamma(y)
=\int _{{\mathbb R}^{n} }\mathbf{M}[ p(t;x/s;y(s;z))\rho (0;z)]d\Gamma(z),}
\end{array}
$$
где $y(s;z)$ -- решение обобщенного уравнения Ито \eqref{Ayd01.1}.

Поскольку левая часть равенства не зависит от $s$, то после дифференцирования по $s$, используя обобщенную формулу Ито, приходим к равенству
$$
 0=\int _{{\mathbb R}^{n} }\mathbf{M}[ d_{s} p(t;x/s;y(s;z))\rho (0;z)]d\Gamma(z)=
$$
$$
=\int _{{\mathbb R}^{n} }\mathbf{M}\Bigl\{\Bigl[ \frac{\partial }{\partial s} p(t;x/s;y(s;z))+a_{j} (s;y(s;z))\frac{\partial }{\partial y_{j} } p(t;x/s;y(s;z))+
$$
$$
\displaystyle{+\frac{1}{2} b_{j,k} (s;y(s;z))b_{ik} (s;y(s;z))\frac{\partial ^{2} }{\partial y_{j} \partial y_{i} } p(t;x/s;y(s;z))\Bigr]ds+}
$$
$$
+\displaystyle b_{j,k} (s;y(s;z))\frac{\partial }{\partial y_{j} } p(t;x/s;y(s;z))dw_{k} (s)+
$$
$$
+\displaystyle\int [p\bigl(t;x/s;y(s;z)+g(t;y(s;z);\gamma )\bigr)- p(t;x/s;y(s;z))]\nu (ds;d\gamma )\Bigr\} \rho (0;z)d\Gamma(z).
$$

Учитывая вновь свойство \eqref{ya}, переходим к следующему равенству:
\begin{equation*}
\begin{array}{c}
\displaystyle\int _{{\mathbb R}^{n} }\Bigl\{ \Bigl[ \frac{\partial }{\partial s} p(t;x/s;y)+a_{j} (s;y)\frac{\partial }{\partial y_{j} } p(t;x/s;y)+
\\
+\displaystyle\frac{1}{2} b_{j,k} (s;y)b_{ik} (s;y)\frac{\partial ^{2} }{\partial y_{j} \partial y_{i} } p(t;x/s;y)\Bigr]
+
\\
+\displaystyle\int \Bigl[p(t;x/s;y+g(s;y;\gamma ))-p(t;x/s;y)\Bigr]  \Pi (d\gamma )\Bigr\} p(s;y)d\Gamma(y)=0.
\end{array}
 \end{equation*}

В свою очередь, требование независимости $p(t;x/s;y)$ от $p(s;y)$  будет выполнено, если
$$
\begin{array}{c}
  \displaystyle \frac{\partial }{\partial s} p(t;x/s;y)+a_{j} (s;y)\frac{\partial }{\partial y_{j} } p(t;x/s;y)+ \\
  +\displaystyle \frac{1}{2} b_{j,k} (s;y)b_{ik} (s;y)\frac{\partial ^{2} }{\partial y_{j} \partial y_{i} } p(t;x/s;y) + \\
  +\displaystyle \int \Bigl[p(t;x/s;y+g(s;y;\gamma ))-p(t;x/s;y)\Bigr]  \Pi (d\gamma )=0.
\end{array}
$$

Полученное уравнение соответствует {\it обратному уравнению Колмогорова } для плотностей переходных вероятностей.

Возвратимся к уравнениям \eqref{Ayd01.1}. Заметим, что при построении уравнений этого параграфа привлекались обобщенные уравнения  Ито \eqref{Ayd01.1} со случайной пуассоновской мерой. Для того  чтобы иметь возможность воспользоваться полученными выше выводами для решений обобщенных уравнений Ито с центрированными мерами, в исходных уравнениях \eqref{Ayd01.1}  потребуется перейти от таких уравнений к уравнениям с пуассоновскими мерами. Это приведет в полученных выше выражениях к замене коэффициентов $a_{j} (\tau ;z)$ на
$$
\bar{a}_{j} (\tau ;z)=a_{j} (\tau ;z)-\displaystyle\int g(\tau ;z;\gamma ) \Pi (d\gamma )
$$
и совпадению с видом уравнений,  полученных И. И. Гихманом и А.~В.~Скороходом \cite[c.~301--302]{GS_68}.

\section{Заключение}

Введение понятия стохастического ядра стохастического интегрального инварианта позволяет получить известные результаты (например, формула Ито-Вентцеля \cite{D_89}, уравнения Колмогорова) что подтверждает полноту и замкнутость теории интегральных инвариантов. Кроме того, метод интегральных инвариантов, рассмотренный с статье, дает возможность дальнейшего развития теории стохастических дифференциальных уравнений и ее применения (например, первые и стохастические первые интегралы, программное управление с вероятностью единица \cite{D_78,D_98,12_KchDubPrep,09_ChUprW,11_KchUpr}).


 \setlength\parsep{0pt} 


\end{document}